\documentclass{article}
\usepackage[a4paper]{geometry}
\usepackage{amsmath,amssymb,amsthm}
\usepackage{mathtools}
\usepackage{tikz-cd}
\usepackage[initials]{amsrefs}
\usepackage{hyperref}
\hypersetup{colorlinks}
\usepackage[english]{babel}

\newtheorem{theorem}{Theorem}[section]
\newtheorem{definitie}[theorem]{Definition} 
\newtheorem{example}[theorem]{Example}
\newtheorem{cor}[theorem]{Corollary}
\newtheorem{prop}[theorem]{Proposition}
\newtheorem{lemma}[theorem]{Lemma}
\newtheorem{rem}[theorem]{Remark}

\newcommand{\bb}[1]{\mathbb{#1}}
\newcommand{\id}{\textup{id}}
\newcommand{\im}{\textup{im}}
\newcommand{\Ima}{\textup{Im}}
\newcommand{\Rea}{\textup{Re}}
\newcommand{\sh}{\textup{Sh}}
\newcommand{\supg}{^{\mathfrak{g}}}
\newcommand{\subg}{_{\mathfrak{g}}}
\newcommand{\supv}{\scriptstyle ^V \displaystyle}
\newcommand{\supu}{\scriptstyle ^U \displaystyle}
\newcommand{\supe}[1]{\scriptstyle ^{#1} \displaystyle}
\newcommand{\osum}[1]{\smashoperator{\sum_{#1}}}
\newenvironment{myproof}{\proof[\bfseries \textup{Proof.\ }]}{\endproof}

\begin{document}
	\title{\vspace{-2.0cm}\textbf{Classification of 2-term $L_\infty$-algebras}}
	\author{\textbf{Kevin van Helden}\footnote{	email: k.s.van.helden@rug.nl}}
	\maketitle
	\vspace{-0.7cm}
	\begin{center}
	Bernoulli Institute for Mathematics, Computer Science and Artificial Intelligence\\
	University of Groningen\\
	Nijenborgh 9,
	9747 AG Groningen,
	The Netherlands
\end{center}
\begin{center}
	Van Swinderen Institute for Particle Physics and Gravity\\\
	University of Groningen\\
	Nijenborgh 4,
	9747 AG Groningen,
	The Netherlands\\
\end{center}

\begin{abstract}
 We classify all 2-term $L_\infty$-algebras up to isomorphism. We show that such $L_\infty$-algebras are classified by a Lie algebra, a vector space, a representation (all up to isomorphism) and a cohomology class of the corresponding Lie algebra cohomology.
\end{abstract}

{
	\hypersetup{linkcolor=black}
	\tableofcontents
}

\section{Introduction}
$L_\infty$-algebras are generalizations of Lie algebras. They are chain complexes of vector spaces on which there is a graded anti-symmetric bracket which satisfies Jacobi-like identities. They were introduced in the early 90's and their applications have become numerous over the last thirty years \cite{Stasheff2018}. In theoretical physics, $L_\infty$-algebras are often an appropriate underlying structure to encode the gauge parameters and the field equations of a field theory \cite{field_theories}. Particular examples of such field theories include closed bosonic string theories, higher Chern-Simons theories and local prequantum field theories \cites{zwiebach,Roytenberg2002,baezhoffrog}.

In classical string theory, another similar structure can also be considered: a Lie 2-algebra. This structure categorifies that of a Lie algebra \cite{Baez2004}. This generalization of a Lie algebra is a natural consequence of replacing point particles in classical mechanics by strings. Whereas point particles can be represented canonically by an object in a category, it is natural to represent strings by a morphism in another suitable 2-category \cite{baezhoffrog}. As both an $L_\infty$-algebra and a Lie 2-algebra describe classical string theory, it is natural to assume that there is a connection between them. This connection was proven by Baez and Crans in 2004. They showed in \cite{Baez2004} that the categories of semistrict Lie 2-algebras and of 2-term $L_\infty$-algebras are equivalent. Moreover, they also proved some primary results into the classification of Lie 2-algebras. One of those results is \cite{Baez2004}*{Theorem 57}, in which they state that there is a one-to-one correspondence between equivalence classes of Lie 2-algebras and isomorphism classes of triples consisting of a Lie algebra $\mathfrak{g}$, a representation $(\rho,V)$ of $\mathfrak{g}$, and a 3-cocycle in the Lie algebra cohomology on $\mathfrak{g}$ with values in $V$.

They do not use this result to classify all semistrict Lie 2-algebras up to isomorphism, however. More recent research has been more focused on string Lie 2-algebras and on applying this structure to string theory than on further classification \cites{Schmidt2019,Ritter2014,Liu2014}. Consequently, further classification of semistrict Lie 2-algebras has not been done yet.

The previous result does however indicate what classification up to isomorphism might look like. As an equivalence of semistrict Lie 2-algebras implies that the two considered structures are isomorphic on homology, the vector spaces that are modded out by taking homology cannot be retrieved from the previously mentioned isomorphism classes of triples. A natural suggestion would thus be to include those in the classification data. By the first isomorphism theorem, those vector spaces must be isomorphic. This would thus require us to extend the isomorphism classes of triples to those of quadruples.

In this paper, we show that this suffices to classify all 2-term $L_\infty$-algebras up to isomorphism. As it was proven in \cite{Baez2004} that the categories of 2-term $L_\infty$-algebras and semistrict Lie 2-algebras are equivalent, our classification of 2-term $L_\infty$-algebras will also classify all semistrict Lie 2-algebras. We choose to prove our statements in terms of 2-term $L_\infty$-algebras, as they allow for simpler and more explicit calculations. We show that the isomorphism classes of all 2-term $L_\infty$-algebras are given by isomorphic quadruples of a Lie algebra, a vector space, a representation, and a cohomology class of the given Lie algebra cohomology.

In Section \ref{secdef}, we will recall the category of 2-term $L_\infty$-algebras and highlight the involved equations, as well as giving some examples of objects and morphisms. In Section \ref{secclass}, we will prove our classification statements in three steps. Firstly, we will introduce elementary concepts from Lie algebra cohomology and then we give a type of example of 2-term $L_\infty$-algebras. Secondly, we show that every 2-term $L_\infty$-algebra is isomorphic to a 2-term $L_\infty$-algebra of the type considered in the example and, thirdly, that those examples are unique up to isomorphic quadruples of a Lie algebra $\mathfrak{g}$, a vector space $U$, a representation $(\rho,V)$, and an element in the Lie algebra cohomology $H^3(\mathfrak{g},\rho,V)$. We will also show how this result relates to \cite{Baez2004}*{Theorem 57}. We conclude by discussing some possible paths for future classification.

\section{The category of 2-term \texorpdfstring{$L_\infty$}{L∞}-algebras}
\label{secdef}
A 2-term $L_\infty$-algebra is an $L_\infty$-algebra concentrated in degree one and zero. We recall the notions of a  2-term $L_\infty$-algebra and a morphism between two of them, and we will give examples for both. We choose our conventions so that our definitions of coincide with (suitable restrictions of those in) \cite{Lada1995} and \cite{Lada1993}. If not explicitly mentioned, we assume all maps between vector spaces are linear.

\begin{definitie}
	\label{lowdimsemistrictdeflie2alg}
	A \textbf{2-term $L_\infty$-algebra} is a graded vector space $L=L_0\oplus L_1$ with 
	\begin{enumerate}
		\item a differential $d\colon L_1 \to L_0$;
		\item a bracket $[\cdot,\cdot]\colon L\wedge L \to L$ of degree 0 such that
		\begin{align}
			\label{lowdimsemistrictdeflie2algeq1}
			d([x,v]) = [x,d(v)]
		\end{align}
	and
	\begin{align}
		\label{lowdimsemistrictdeflie2algeq1b}
	[d(u),v] = [u,d(v)]
	\end{align}
		for $x\in L_0$ and $u,v\in L_1$;
		\item a Jacobiator $J\colon L_0\wedge L_0\wedge L_0\to L_1 $ such that
		\begin{align}
			\label{lowdimsemistrictdeflie2algeq2}
			d(J_{x,y,z}) = [x,[y,z]]- [[x,y],z]-[y,[x,z]]
		\end{align}
		and
		\begin{align}
			\label{lowdimsemistrictdeflie2algeq3}
			J_{d(v),y,z} = [v,[y,z]]- [[v,y],z]-[y,[v,z]]
		\end{align}
		for $x,y,z\in L_0$ and for $v\in L_1$ (note that we write $J_{x,y,z}:=J(x\wedge y\wedge z)$),
	\end{enumerate}
	such that the following equation holds:
\begin{align}
\label{lowdimsemistrictdeflie2algeq4}
\osum{\sigma \in\sh(1,3)}(-1)^\sigma[x_{\sigma(1)},J_{x_{\sigma(2)},x_{\sigma(3)},x_{\sigma(4)}}]\
-\osum{\sigma\in \sh(2,2)}(-1)^\sigma J_{[x_{\sigma(1)},x_{\sigma(2)}],x_{\sigma(3)},x_{\sigma(4)}}=0
\end{align}
for $x_1,x_2,x_3,x_4 \in L_0$, where 
\begin{align*}
\sh(m,n):= \{ \sigma \in S_{m+n} \ | \ \sigma(i)<\sigma(i+1) \ \textup{for all} \ 1\leq i \leq m+n-1, i\neq m\}
\end{align*}
for $m,n\in \bb{N}$.
\end{definitie}
It is worthwhile to note that the bracket and the Jacobiator can also be seen as a chain map and a chain homotopy on their respective chain complexes. The interested reader can find more about this in \cite{Roytenberg2007}. Furthermore, we note for the reader familiar with \cites{Lada1995,Lada1993} that the sign involved in Equation \eqref{lowdimsemistrictdeflie2algeq4} is the ordinary permutation sign, in contrast to the Koszul sign used in \cites{Lada1995,Lada1993}. This is due to the fact that the other signs that are absorbed into the Koszul sign are identically one in the 2-term $L_\infty$-algebra case. Normally, these additional signs would appear as a consequence of permuting multiple odd degree elements, but, by degree considerations, such terms are already zero. 
\begin{example}
	\label{examplelie2alg}
	\textup{ A simple example of a 2-term $L_\infty$-algebra can be given as follows. Let $\bb{H}$ be the algebra of quaternions, and for $z\in \bb{H}$, let $\Rea(z)$ and $\Ima(z)$ denote the real and imaginary part of $z$ respectively, that is, $\Rea(1)=1$, $\Ima(\mathbf{i})=\mathbf{i}$, $\Ima(\mathbf{j})=\mathbf{j}$, $\Ima(\mathbf{k})=\mathbf{k}$ and $\Ima(1)=\Rea(\mathbf{i})=\Rea(\mathbf{j})=\Rea(\mathbf{k})=0$. Then for $v\in \bb{H}$, the graded vector space $\bb{H}\oplus \bb{H}$ is a 2-term $L_\infty$-algebra, with
		\begin{enumerate}
			\item a differential
			$	d=\Rea\colon \bb{H} \to \bb{H};$
			\item a bracket $[\cdot,\cdot]\colon \bb{H} \oplus \bb{H} \wedge \bb{H} \oplus \bb{H} \to \bb{H} \oplus \bb{H}$ given by
			\begin{align*}
			[a\oplus b,c\oplus d] = \Ima(\Ima(a)\Ima(c))\oplus\Ima(\Ima(a)\Ima(d)+\Ima(b)\Ima(c));
			\end{align*}
			\item 	a Jacobiator $J\colon \Lambda^3 \bb{H} \to \bb{H}$ determined by
			\begin{center}
			\begin{tabular}{c c c c}
			$J(1,\mathbf{i},\mathbf{j})=0;$ & $J(1,\mathbf{i},\mathbf{k}) =0;$ & $J(1,\mathbf{j},\mathbf{k})=0;$ & $J(\mathbf{i},\mathbf{j},\mathbf{k})=\Ima(v).$
			\end{tabular}
			\end{center}
		\end{enumerate}
	}
\end{example}
The next two examples come from \cite{Baez2004} and \cite{Baez2007} respectively.
\begin{example}
	\label{exampleskeletalstringlie2alg}
	\textup{ Let $\mathfrak{g}$ be a Lie algebra with Lie bracket $[\cdot,\cdot]_{\mathfrak{g}}$, and let $k\in \bb{R}$. The \textbf{skeletal string 2-term $L_\infty$-algebra} consists of a graded vector space $\mathfrak{g}_k:=\mathfrak{g}\oplus \bb{R}$ with
		\begin{enumerate}
			\item a differential $d=0\colon\bb{R}\to \mathfrak{g}$;
			\item a bracket $[\cdot,\cdot]\colon \mathfrak{g}_k \wedge \mathfrak{g}_k \to \mathfrak{g}_k$ given by 
			\begin{align*}
				[x,y] = [x,y]_{\mathfrak{g}}
			\end{align*} and 
		\begin{align*}
			[x,v]=0
		\end{align*}
	 for $x,y\in\mathfrak{g}$ and $v\in \bb{R}$;
			\item a Jacobiator $J\colon\Lambda^3 \mathfrak{g}_k \to \bb{R}$ given by
			\begin{align*}
				J(x,y,z) = \langle x,[y,z] \rangle,
			\end{align*}
		where $\langle \cdot,\cdot \rangle$ is the Killing form on $\mathfrak{g}$.
		\end{enumerate}
}
\end{example}
\begin{example}
	\label{examplestrictstringlie2alg}
	\textup{ Let $\mathfrak{g}$ be a Lie algebra with Lie bracket $[\cdot,\cdot]_{\mathfrak{g}}$, and let $k\in \bb{R}$. Furthermore, define the (real) vector spaces
		\begin{align*}
			\mathcal{P}_0\mathfrak{g} &= \{f\in C^\infty([0,2\pi],\mathfrak{g}) \ | \ f(0)=0\},\\
			\Omega\mathfrak{g} &= \{f\in C^\infty(\bb{R}/2\pi\bb{Z},\mathfrak{g}) \ | \ f(0)=0\},
		\end{align*}
		and the quotient map $q\colon [0,2\pi] \to \bb{R}/2\pi\bb{Z}$. Moreover, note that $\mathcal{P}_0\mathfrak{g}$ is a Lie algebra with the Lie bracket given by 
		\begin{align*}
			[f,g]_{	\mathcal{P}_0\mathfrak{g}} &:= (\theta\mapsto  [f(\theta),g(\theta)]_{\mathfrak{g}}).
		\end{align*}
		The \textbf{strict string 2-term $L_\infty$-algebra} consists of a graded vector space $\mathcal{P}_k\mathfrak{g}:=\mathcal{P}_0\mathfrak{g}\oplus (\Omega\mathfrak{g} \oplus \bb{R})$ with
		\begin{enumerate}
			\item a differential 
			\begin{align*}
				\begin{split}
					d\colon \Omega\mathfrak{g} \oplus \bb{R} &\to \mathcal{P}_0\mathfrak{g}\\
					f \oplus k &\mapsto q^*(f)
				\end{split}
			\end{align*}
			\item a bracket $[\cdot,\cdot]\colon \mathcal{P}_k\mathfrak{g} \wedge \mathcal{P}_k\mathfrak{g} \to \mathcal{P}_k\mathfrak{g}$ given by 
			\begin{align*}
				[f,g] = [f,g]_{	\mathcal{P}_0\mathfrak{g}}
			\end{align*} and 
			\begin{align*}
				[f,h\oplus l]=[f,q^*(h)]_{	\mathcal{P}_0\mathfrak{g}} \oplus 2k\int_0^{2\pi} \langle f(\theta),h'(\theta)\rangle d\theta 
			\end{align*}
			for $f,g\in	\mathcal{P}_0\mathfrak{g}$ and $h\oplus l\in \Omega\mathfrak{g} \oplus \bb{R}$, where $\langle \cdot,\cdot \rangle$ is the Killing form on $\mathfrak{g}$;
			\item a Jacobiator $J=0\colon\Lambda^3 \mathcal{P}_0\mathfrak{g} \to \Omega\mathfrak{g} \oplus \bb{R}$.
		\end{enumerate}
	}
\end{example}
It is proven in \cite{Baez2007} that $\mathfrak{g}_k$ and $\mathcal{P}_k\mathfrak{g}$ are equivalent as 2-term $L_\infty$-algebras. It is not hard to show that both 2-term $L_\infty$-algebras are not isomorphic, but, for expository purposes, we postpone this to the end of the paper.
\begin{definitie}
	\label{lowdimsemistrictlie2aglmor}
	Let $(L,d,[\cdot,\cdot],J)$ and $(L',d',[\cdot,\cdot]',J')$ be two 2-term $L_\infty$-algebras. A 2-term \textbf{$L_\infty$-algebra morphism} from $L$ to $L'$ is a pair $(\phi,\Phi)$, with
	\begin{itemize}
		\item a map $\phi\colon L \to L '$ of degree 0 such that
		\begin{align}
			\label{lowdimsemistrictlie2aglmoreq1}
			\phi(d(v)) = d'(\phi(v))
		\end{align}
		for $v\in L_1$;
		\item a map $\Phi\colon L_0\wedge L_0 \to L_1'$ such that
		\begin{align}
			\label{lowdimsemistrictlie2aglmoreq2}
			d'(\Phi(x\wedge y)) = \phi([x,y])-[\phi(x),\phi(y)]'
		\end{align}
		and
		\begin{align}
			\label{lowdimsemistrictlie2aglmoreq3}
			\Phi(d(v)\wedge y) = \phi([v,y]) - [\phi(v),\phi(y)]'
		\end{align}
		for $x,y\in L_0$ and $v\in L_1$,
	\end{itemize} 
	such that the following equation holds:
	\begin{align}
\label{lowdimsemistrictlie2algmoreq4}
\phi(J_{x_1,x_2,x_3})-J'_{\phi(x_1),\phi(x_2),\phi(x_3)}
= \osum{\sigma \in \sh(1,2)}(-1)^\sigma\big([\phi(x_{\sigma(1)}),\Phi_{x_{\sigma(2)},x_{\sigma(3)}}]'+\Phi_{x_{\sigma(1)},[x_{\sigma(2)},x_{\sigma(3)}]}\big)
\end{align}
for $x_1,x_2,x_3 \in L_0$, where $\Phi_{x_1,x_2} := \Phi(x_1 \wedge x_2)$.
\end{definitie}
Again, the maps $\phi$ and $\Phi$ are a chain map and chain homotopy on their respective chain complexes. More information on this can be found in \cite{Roytenberg2002}. Similar as in Equation \eqref{lowdimsemistrictdeflie2algeq4}, the Koszul sign the authors of \cites{Lada1995,Lada1993} used in the general version of Equation \eqref{lowdimsemistrictlie2algmoreq4} reduces to a permutation sign due to the lack of non-trivial elements of degree higher than one.
 
As a morphism between 2-term $L_\infty$-algebras consists of two different maps, it might not be evident how to compose them. We will thus give the description of composition of said morphisms in the following definition.
\begin{definitie}
	\label{lowdimsemistrictlemchahom}
	Let $(L,d,[\cdot,\cdot],J)$, $(L',d',[\cdot,\cdot]',J')$ and $(L'',d'',[\cdot,\cdot]'',J'')$ be three 2-term $L_\infty$-algebras, and let 
	\[
	(\phi,\Phi)\colon (L,d,[\cdot,\cdot],J)\to (L',d',[\cdot,\cdot]',J')\]
	and 
	\[
	(\phi',\Phi')\colon (L',d',[\cdot,\cdot]',J')\to (L'',d'',[\cdot,\cdot]'',J'')
	\]
	 be two 2-term $L_\infty$-algebra morphisms. The \textbf{composition} of $(\phi,\Phi)$ and $(\phi',\Phi')$ is $(\phi'\circ \phi, \Psi_{
	\Phi',\Phi})$, where $\Psi_{\Phi',\Phi}$ is the map  
	\begin{align*}
		\begin{split}
		\Psi_{\Phi',\Phi}\colon L_0 \wedge L_0 &\to L_1''\\
		x\wedge y &\mapsto \Phi'_{\phi(x),\phi(y)} +\phi'(\Phi_{x,y}).
		\end{split}
	\end{align*} 
\end{definitie}
It is clear from this definition that the identity morphism is given by $(\id,0)$. It is shown in \cite{Liu2014} that a morphism of 2-term $L_\infty$-algebras $(\phi,\Phi)$ is an isomorphism of 2-term $L_\infty$-algebras if $\phi$ is a linear isomorphism, as the inverse morphism is given by $(\phi^{-1},\Phi')$, where $\Phi'\colon L_0'\wedge L_0'\to L_1$ is given by $\Phi'(x,y)= -\phi^{-1}\Phi(\phi^{-1}(x),\phi^{-1}(y))$.
\begin{example}
	\label{examplelie2algmor}
	\textup{ Consider the 2-term $L_\infty$-algebra $\bb{H}\oplus \bb{H}$ from Example \ref{examplelie2alg}. Then the two maps
		\begin{enumerate}
			\item $\phi\colon \bb{H}\oplus \bb{H} \to \bb{H}\oplus \bb{H}$ determined by
			\begin{center}
				\begin{tabular}{ c c c c }
					$\phi(1\oplus 0)=1\oplus 0;$ & $\phi(\mathbf{i}\oplus 0)=\mathbf{j}\oplus 0;$ & $\phi(\mathbf{j}\oplus 0)=\mathbf{k}\oplus 0;$ & $\phi(\mathbf{k}\oplus 0)=\mathbf{i}\oplus 0;$\\ 
				$\phi(0\oplus 1)=0\oplus 1$ & $\phi(0\oplus \mathbf{i})=0\oplus \mathbf{j};$ & $\phi(0\oplus \mathbf{j})=0\oplus \mathbf{k};$ & $\phi(0\oplus \mathbf{k})=0\oplus \mathbf{i}$,
				\end{tabular}
			and
			\end{center}
			\item $\Phi\colon \Lambda^2 \bb{H} \to \bb{H}$ determined by
			\begin{center}
				\begin{tabular}{c c c}
					$\Phi(1\wedge \mathbf{i})=0;$ & $\Phi(1\wedge \mathbf{k})=0;$ & $\Phi(\mathbf{j}\wedge \mathbf{k})=\Rea(v(\mathbf{j}-\mathbf{k}))\mathbf{i};$\\
					$\Phi(1\wedge \mathbf{j})=0;$ & $\Phi(\mathbf{i}\wedge \mathbf{j})=\Rea(v(\mathbf{i}-\mathbf{j}))\mathbf{k};$ & $\Phi(\mathbf{k}\wedge \mathbf{i})=\Rea(v(\mathbf{k}-\mathbf{i}))\mathbf{j}$
				\end{tabular}
			\end{center}
			form an automorphism of 2-term $L_\infty$-algebras on $\bb{H}\oplus \bb{H}$.
		\end{enumerate}
	}
\end{example}

\section{Classification}
\label{secclass}
In this section, we will give the full classification procedure of 2-term $L_\infty$-algebras. For the sake of clarity, this section is split into two subsections.

\subsection{The fundamental example}
In this subsection, we will define some important aspects of Lie algebra cohomology as in \cite{Baez2004} and we give an example of a 2-term $L_\infty$-algebra using this cohomology. This type of example will be quintessential in the two proofs in the next subsection.

Let $\mathfrak{g}$ be a Lie algebra with Lie bracket $[\cdot,\cdot]_\mathfrak{g}$, and let $(\rho,V)$ be a Lie algebra representation of $\mathfrak{g}$. We then define the cochain complex $C(\mathfrak{g},\rho,V)$ by setting $C_n(\mathfrak{g},\rho,V)$ to be the vector space of linear maps $f\colon\Lambda^n \mathfrak{g} \to V$ and by setting the differential $\delta$ to be 
\begin{multline}
	(\delta f)(x_1,\dots,x_{n+1}):= \osum{\sigma \in \sh(1,n)} (-1)^\sigma\rho(x_{\sigma(1)})f(x_{\sigma(2)},\dots,x_{\sigma(n+1)})\\
	 -\osum{\sigma\in \sh(2,n-1)} (-1)^\sigma f([x_{\sigma(1)},x_{\sigma(2)}]_\mathfrak{g},x_{\sigma(3)},\dots,x_{\sigma(n+1)})
\end{multline}
for $f\in C_n(\mathfrak{g},\rho,V)$. It has been already proven in \cite{Chevalley1948} that $\delta^2=0$. We call the cohomology of this complex $H(\mathfrak{g},\rho,V)$.
\begin{definitie}
	\label{defncocycle}
	An \textbf{$n$-cocycle} is an element $f\in C_n(\mathfrak{g},\rho,V)$ such that $\delta f =0$.
\end{definitie}

\begin{definitie}
	Let $\mathfrak{g}$ and $\mathfrak{h}$ be two Lie algebras and let $(\rho,V)$ and $(\sigma,W)$ be a Lie algebra representation of $\mathfrak{g}$ and of $\mathfrak{h}$ respectively. We say that $J\in C_n(\mathfrak{g},\rho,V)$ and $K\in  C_n(\mathfrak{h},\sigma,W)$ are \textbf{cohomologous} if there exists a Lie algebra morphism $\psi\colon \mathfrak{g}\to \mathfrak{h}$ and an intertwiner $t\colon V \to W$ and a linear map $\Phi\colon \Lambda^{n-1} \mathfrak{g} \to W$ such that 
	\begin{align}
		\label{eqcohomologous3cocycles}
		t(J(x_1,\dots,x_n)) - K(\psi(x_1),\dots,\psi(x_n)) =(\delta \Phi)(x_1,\dots,x_n)
	\end{align}
for all $x_1,\dots,x_n\in \mathfrak{g}$, where $\delta$ in Equation \eqref{eqcohomologous3cocycles} belongs to the cochain complex $C(\mathfrak{g},\sigma \circ \psi,W)$.
\end{definitie}

Note that stating that $J$ and $K$ are cohomologous is equivalent to the stating that $t\circ J$ and $\psi^*K$ are representatives of the same cohomology class in $H(\mathfrak{g},\sigma\circ \psi,W)$.

Using those definitions from Lie algebra cohomology, we give a type of example of 2-term $L_\infty$-algebra, which will play a crucial role in the further classification. In the following lemma and the remainder of the article, we will use the following notation for the sake of clarity. If $U$ and $V$ are vector spaces, and $x\in U \oplus V$, we denote the projection of $x$ on $U$ by $x\supu$.

\begin{lemma}
	\label{lowdimlemexamplelie2alg}
	Let $\mathfrak{g}$ be a Lie algebra with Lie bracket $[\cdot,\cdot]_\mathfrak{g}$, let $U$ be a vector space, let $(\rho,V)$ be a Lie algebra representation of $\mathfrak{g}$, and let $\tilde{J}$ be a 3-cocycle of $(\rho,V)$ on  $\mathfrak{g}$. We then define the graded vector space $L=L^{\mathfrak{g},U,\rho,\tilde{J}}$ by setting $L_0 := \mathfrak{g}\oplus U$ and $L_1:= V\oplus U$. We furthermore define the linear maps
	\begin{enumerate}
		\item $d\colon L_1 \to L_0$ given by $d(v) =v\supu$;
		\item $[\cdot,\cdot]\colon L\wedge L\to L$ given by $[x,y] = [x\supg,y\supg]\subg $, by $[x,v] = -[v,x] = \rho(x\supg)v\supv$ and by $[u,v]=0$;
		\item $J\colon L_0\wedge L_0\wedge L_0 \to L_1$ by $J_{x,y,z}:= \tilde{J}_{x\supg,y\supg,z\supg}$
	\end{enumerate}
	for  $x,y,z \in L_0$ and $u,v \in L_1$. Then $(L,d,[\cdot,\cdot],J)$ is a 2-term $L_\infty$-algebra.
\end{lemma}
\begin{myproof} 
Let $x,y,z,w\in L_0$ and $u,v\in L_1$. Some simple observations, such as the fact that $d(V)=0$ and that $d(L_1)\supg =0$, yield that
\begin{align*}
	d([x,v]) = d\left(\rho(x\supg)v\supv\right) = 0 = [x\supg, d(v)\supg]_\mathfrak{g} = [x, d(v)],
\end{align*}
and that
\begin{align*}
	\begin{split}
	[d(u),v] =  \rho(d(u)\supg)v\supv=0
	=-\rho(d(v)\supg)u\supv
	=[u, d(v)].
	\end{split}
\end{align*}
By furthermore noting that $[x,y] = [x,y]\supg$ and using the Jacobi identity on $\mathfrak{g}$, we also find that
\begin{align*}
	\begin{split}		
	d(J_{x,y,z}) &= d(\tilde{J}(x\supg,y\supg,z\supg))
	= 0 
	= [x\supg,[y\supg,z\supg]_\mathfrak{g}]_\mathfrak{g}-[[x\supg,y\supg]_\mathfrak{g},z\supg]_\mathfrak{g} - [y\supg,[x\supg,z\supg]_\mathfrak{g}]_\mathfrak{g}\\
	&= [x\supg,[y,z]\supg]_\mathfrak{g}-[[x,y]\supg,z\supg]_\mathfrak{g} - [y\supg,[x,z]\supg]_\mathfrak{g}
	= [x,[y,z]]-[[x,y],z] - [y,[x,z]].
	\end{split}
\end{align*}
Similarly, by also using that $(\rho,V)$ is a representation on $\mathfrak{g}$ and that $[x,v]=[x,v]\supv$, we obtain
\begin{align*}
	\begin{split}
		J_{d(v),y, z} &= \tilde{J}(d(v)\supg,y\supg,z\supg)= 0 = -\rho([y\supg,z\supg]_\mathfrak{g})v\supv -\rho (z\supg)\rho(y\supg)v\supv +\rho(y\supg)\rho(z\supg)v\supv\\
	&= -\rho([y,z]\supg)v\supv +\rho (z\supg)[v,y]\supv  - \rho(y\supg)[v,z]\supv
	= [v,[y,z]]-[[v,y],z] - [y,[v,z]].
	\end{split}
\end{align*}
Moreover, we find that
\begin{align*}
	[x,J_{y,z,w}] = \rho(x\supg)J_{y,z,w}\supv= \rho(x\supg)\tilde{J}_{y\supg,z\supg,w\supg}
\end{align*}
and that
\begin{align*}
	J_{[x,y],z,w} = \tilde{J}_{[x,y]\supg,z\supg,w\supg} = \tilde{J}_{[x\supg,y\supg]_\mathfrak{g},z\supg,w\supg},
\end{align*}
giving us for $x_1,x_2,x_3,x_4\in L_0$ that
\begin{align*}
\begin{split}
		\osum{\sigma \in\sh(1,3)}&(-1)^\sigma[x_{\sigma(1)},J_{x_{\sigma(2)},x_{\sigma(3)},x_{\sigma(4)}}]\
		-\osum{\sigma\in \sh(2,2)}(-1)^\sigma J_{[x_{\sigma(1)},x_{\sigma(2)}],x_{\sigma(3)},x_{\sigma(4)}}\\
		&= \osum{\sigma \in\sh(1,3)}(-1)^\sigma\rho(x_{\sigma(1)}\supg)\tilde{J}_{x_{\sigma(2)}\supg,x_{\sigma(3)}\supg,x_{\sigma(4)}\supg}  -\osum{\sigma\in \sh(2,2)}(-1)^\sigma \tilde{J}_{[x_{\sigma(1)}\supg,x_{\sigma(2)}\supg]_\mathfrak{g},x_{\sigma(3)}\supg,x_{\sigma(4)}\supg}\\
		&= \delta \tilde{J}(x_1\supg,x_2\supg,x_3\supg,x_4\supg),
\end{split}
\end{align*}
which is equal to zero as $\tilde{J}$ is a 3-cocycle. This proves that Equations \eqref{lowdimsemistrictdeflie2algeq1}-\eqref{lowdimsemistrictdeflie2algeq4} hold.
\end{myproof}
In this construction, the vector space $U$ appears twice: once in degree 0, and once in degree 1. It should be clear from the context in which degree elements of $U$ reside. If this is not the case, we will explicitly mention the relevant degrees.
\begin{example}
	\textup{
Example \ref{examplelie2alg} can be viewed as a 2-term $L_\infty$-algebra of the shape of Lemma \ref{lowdimlemexamplelie2alg}. First, we note that the imaginary part of $\bb{H}$, $\Ima(\bb{H})$, is a Lie algebra with bracket $[a,b]:=\Ima(\Ima(a)\Ima(b))$, which is isomorphic to $\mathfrak{so}(3)$. By setting $\mathfrak{g}$ equal to $\Ima(\bb{H})$, by setting $U$ equal to the real part of $\bb{H}$, by setting $(\rho,V)$ equal to the adjoint representation of $\Ima(\bb{H})$, and by setting $\tilde{J}\colon\Lambda^3\Ima(\bb{H}) \to \Ima(\bb{H})$ to be the map that sends $\mathbf{i}\wedge \mathbf{j} \wedge \mathbf{k}$ to $\Ima(v)\in \Ima(\bb{H})$, we obtain $\bb{H}\oplus \bb{H}$ from Example \ref{examplelie2alg} as $L^{\mathfrak{g},U,\rho,\tilde{J}}$.
}
\end{example}

Now we will show that each 2-term $L_\infty$-algebra gives rise to a 2-term $L_\infty$-algebra in the shape of these described in Lemma \ref{lowdimlemexamplelie2alg}.
\begin{lemma}
	\label{lem2exlie2alg}
	Let  $(L=L_0\oplus L_1,d,[\cdot,\cdot],J)$ be a 2-term $L_\infty$-algebra and fix a vector space decomposition $L_0\oplus L_1= (\mathfrak{g} \oplus \im(d)) \oplus (\ker(d)\oplus U)$.  Then consider the linear isomorphism 
	\begin{align*}
	\begin{split}
	f\colon\ker(d) \oplus U &\to \ker(d)\oplus \im(d)\\
	v\oplus u &\mapsto v\oplus d(u)
	\end{split}
	\end{align*}
	and the map
	\begin{align*}
	\begin{split}
	h\colon L_0 &\to U\\
	x &\mapsto f^{-1}(x\supe{\im(d)}).
	\end{split}
	\end{align*}
	Then the following statements hold:
	\begin{enumerate}
		\item $\mathfrak{g}$ is a Lie algebra with bracket 
		\begin{align*}
			\begin{split}
				[\cdot,\cdot]_\mathfrak{g}\colon \mathfrak{g} \wedge \mathfrak{g} &\to \mathfrak{g}\\
				[y,z]_\mathfrak{g} &= [y,z]\supg;
			\end{split}
		\end{align*} 
	\item the map \begin{align*}
		\begin{split}
			\rho\colon \mathfrak{g} &\to \mathfrak{gl}(\ker(d))\\
			\rho(x)(v) &= [x,v]\supe{\ker(d)}.
		\end{split}
	\end{align*}
	is a Lie algebra representation of $\mathfrak{g}$ on $\ker(d)$;
\item the map
\begin{align*}
	\begin{split}
		\tilde{J}\colon \mathfrak{g} \wedge \mathfrak{g} \wedge \mathfrak{g} &\to \ker(d)\\
		x_1\wedge x_2 \wedge x_3 &\mapsto J_{x_1,x_2,x_3}\supe{\ker(d)}-\sum_{\sigma \in \sh(1,2)}(-1)^\sigma [x_{\sigma(1)},h([x_{\sigma(2)},x_{\sigma(3)}])]^{\ker(d)}
	\end{split}
\end{align*}
is a 3-cocycle of $(\rho,\ker(d))$ on $\mathfrak{g}$.
\end{enumerate}
	\end{lemma}
\begin{myproof}
\begin{enumerate}
	\item Note that the bracket is antisymmetric by definition. Furthermore, we find that
	\begin{align*}
		\begin{split}
			[x,[y,z]_\mathfrak{g} ]_\mathfrak{g} &-[[x,y]_\mathfrak{g},z]_\mathfrak{g} - [y,[x,z]_\mathfrak{g}]_\mathfrak{g}
			=([x,[y,z]\supg]-[[x,y]\supg,z]-[y,[x,z]\supg])\supg\\
			&=([x,[y,z]]-[[x,y],z]-[y,[x,z]])\supg
			=d(J_{x,y,z})\supg=0,
		\end{split}
	\end{align*}
	for $x,y,z\in \mathfrak{g}$, so  $[\cdot,\cdot]_\mathfrak{g}$ satisfies the Jacobi identity and is thus a Lie bracket.  Note that the second equality stems from the fact that, by Equation \eqref{lowdimsemistrictdeflie2algeq1}, $[x,w]\supg = [x,d(u)]\supg = (d([x,u]))\supg=0$ for $u\in L_1$ and $w=d(u)$.
	\item  For $x,y,z\in \mathfrak{g}$ and $v\in \ker(d)$, Equation \eqref{lowdimsemistrictdeflie2algeq1} implies that $d([x,v]) = [x,d(v)]=0$, so $[x,v] \in \ker(d)$. Subsequently, we obtain that
	\begin{align*}
		\begin{split}
			\big(\rho([y,z])&-\rho(y)\rho(z)+\rho(z)\rho(y)\big)(v)
			=\rho([y,z])(v) -\rho(z)([v,y])+\rho(y)([v,z])\\
			&= -[v,[y,z]]+ [[v,y],z]+[y,[v,z]] = -J_{d(v),y,z} =0,
		\end{split}
	\end{align*}
	proving that $(\rho,\ker(d))$ is a representation.
	\item For $u\in U$ and $z\in L_0$, we note that, by applying Equation \eqref{lowdimsemistrictdeflie2algeq1}, we obtain
	\begin{align}
	\label{eq:handu}
	h([x,d(u)]) = h(d([x,u]))=[x,u]\supu = [x,h(d(u))]\supu.
	\end{align}	
	For $x_1,x_2,x_3,x_4\in \mathfrak{g}$, using Equation \eqref{eq:handu}, Equation \eqref{lowdimsemistrictdeflie2algeq2} and Equation \eqref{lowdimsemistrictdeflie2algeq3} yields
	\begin{align*}
	\begin{split}
	&\osum{\sigma \in\sh(1,3)}(-1)^\sigma [x_{\sigma(1)},J_{x_{\sigma(2)},x_{\sigma(3)},x_{\sigma(4)}}\supe{U}] = \osum{\sigma \in\sh(1,3)}(-1)^\sigma [x_{\sigma(1)},h(dJ_{x_{\sigma(2)},x_{\sigma(3)},x_{\sigma(4)}})]\\
	&\ \ = \osum{\sigma \in \sh(1,1,2)}(-1)^\sigma [x_{\sigma(1)},h([x_{\sigma(2)},[x_{\sigma(3)},x_{\sigma(4)}]])] \\
	&\ \ = \osum{\sigma \in \sh(1,1,2)}(-1)^\sigma ([x_{\sigma(1)},h([x_{\sigma(2)},[x_{\sigma(3)},x_{\sigma(4)}]\supg])]+  [x_{\sigma(1)},h([x_{\sigma(2)},[x_{\sigma(3)},x_{\sigma(4)}]\supe{\im(d)}])])\\
	&\ \ = \osum{\sigma \in \sh(1,1,2)}(-1)^\sigma ([x_{\sigma(1)},h([x_{\sigma(2)},[x_{\sigma(3)},x_{\sigma(4)}]\supg])]+  [x_{\sigma(1)},[x_{\sigma(2)},h([x_{\sigma(3)},x_{\sigma(4)}])]\supu])
	\end{split}
	\end{align*}
	and
	\begin{align*}
	\begin{split}
	&\osum{\sigma\in \sh(2,2)}(-1)^\sigma J_{[x_{\sigma(1)},x_{\sigma(2)}]^{\im(d)},x_{\sigma(3)},x_{\sigma(4)}} = \osum{\sigma\in \sh(2,2)}(-1)^\sigma J_{d(h([x_{\sigma(1)},x_{\sigma(2)}])),x_{\sigma(3)},x_{\sigma(4)}}\\
	&\ = 	\osum{\sigma\in \sh(2,2)}(-1)^\sigma [h([x_{\sigma(1)},x_{\sigma(2)}]),[x_{\sigma(3)},x_{\sigma(4)}]]+\osum{\sigma \in \sh(1,1,2)}(-1)^\sigma [x_{\sigma(1)},[x_{\sigma(2)},h([x_{\sigma(3)},x_{\sigma(4)}])]]
	\end{split}
	\end{align*}
	with 
	\begin{align*}
	\sh(1,1,2) = \{ \sigma \in S_4 \ | \ \sigma(3)< \sigma(4)\}.
	\end{align*}
	As a direct consequence of Equation \eqref{lowdimsemistrictdeflie2algeq1}, we find that
	\begin{align*}
	\begin{split}
	[[x_1,x_2]\supe{\im(d)},h([x_3,x_4])]\supe{\ker(d)} = 0.
	\end{split}
	\end{align*}
In turn, we then yield that
\begin{align*}
	\begin{split}
		\delta& \tilde{J}(x_1,x_2,x_3,x_4)\\
		&= \osum{\sigma \in\sh(1,3)}(-1)^\sigma[x_{\sigma(1)},J_{x_{\sigma(2)},x_{\sigma(3)},x_{\sigma(4)}}\supe{\ker(d)}] - \osum{\sigma \in \sh(1,1,2)} (-1)^\sigma[x_{\sigma(1)},[x_{\sigma(2)},h([x_{\sigma(3)},x_{\sigma(4)}])]\supe{\ker(d)}]\\
		&-\osum{\sigma\in \sh(2,2)}(-1)^\sigma J_{[x_{\sigma(1)},x_{\sigma(2)}]^\mathfrak{g},x_{\sigma(3)},x_{\sigma(4)}}
		+ \osum{\sigma\in \sh(2,2)}(-1)^\sigma [[x_{\sigma(1)},x_{\sigma(2)}]\supg,h([x_{\sigma(3)},x_{\sigma(4)}])]\\
		&\ \ +\osum{\sigma \in \sh(1,1,2)} [x_{\sigma(1)},h([x_{\sigma(2)},[x_{\sigma(3)},x_{\sigma(4)}]\supg])]\\
&=	\osum{\sigma \in\sh(1,3)}(-1)^\sigma \big([x_{\sigma(1)},J_{x_{\sigma(2)},x_{\sigma(3)},x_{\sigma(4)}}\supe{\ker(d)}+J_{x_{\sigma(2)},x_{\sigma(3)},x_{\sigma(4)}}\supe{U}]\big)\supe{\ker(d)}\\
&\ \ -\osum{\sigma\in \sh(2,2)}(-1)^\sigma \big(J_{[x_{\sigma(1)},x_{\sigma(2)}]^\mathfrak{g},x_{\sigma(3)},x_{\sigma(4)}}+J_{[x_{\sigma(1)},x_{\sigma(2)}]^{\im(d)},x_{\sigma(3)},x_{\sigma(4)}}\big)\supe{\ker(d)}\\
&=	\bigg(\osum{\sigma \in\sh(1,3)}(-1)^\sigma[x_{\sigma(1)},J_{x_{\sigma(2)},x_{\sigma(3)},x_{\sigma(4)}}]
-\osum{\sigma\in \sh(2,2)}(-1)^\sigma J_{[x_{\sigma(1)},x_{\sigma(2)}],x_{\sigma(3)},x_{\sigma(4)}}\bigg)\supe{\ker(d)},
\end{split}
\end{align*}
which is zero as $J$ satisfies Equation \eqref{lowdimsemistrictdeflie2algeq4}, proving that $\tilde{J}$ is a 3-cocycle of $\mathfrak{g}$ on $(\rho,\ker(d))$.	
\end{enumerate}
\end{myproof}
If we use the newly found Lie algebra, representation and 3-cocycle from Lemma \ref{lem2exlie2alg} in the construction given in Lemma \ref{lowdimlemexamplelie2alg}, we can create a new 2-term $L_\infty$-algebra.
\begin{cor}
	\label{cornewlie2alg}
	Let $L$ be a 2-term $L_\infty$-algebra and fix a vector space decomposition $(\mathfrak{g}\oplus \im(d))\oplus (\ker(d)\oplus U)$. Moreover, let $\mathfrak{g}$ (as a Lie algebra), $(\rho,\ker(d))$ and $\tilde{J}$ as in Lemma \ref{lem2exlie2alg}. Then $L^{\mathfrak{g},\im(d),\rho,\tilde{J}}$ is a 2-term $L_\infty$-algebra. \qed
\end{cor}

\subsection{Main classification theorems}
In this subsection, we will prove the main classification theorems for 2-term $L_\infty$-algebras. In the first theorem, we prove that every 2-term $L_\infty$-algebra is isomorphic to one in the shape of these from Lemma \ref{lowdimlemexamplelie2alg}. We will do this by proving a proposition which states that the structure of a 2-term $L_\infty$-algebra can be transferred to an isomorphic graded vector space, for which the new differential, bracket and Jacobiator are uniquely determined by a given isomorphism of graded vector spaces and a given linear map. This is a special case of the Homotopy Transfer Theorem, which states that the structure of a 2-term $L_\infty$-algebra can be transferred to any of its homotopy retracts \cite{loday2012algebraic}. In general, the 2-term $L_\infty$-algebra structure induced by the Homotopy Transfer Theorem is not isomorphic to the aforementioned structure, but both structures are rather isomorphic on homology. We will focus on the situation in which the homotopy retract is linearly isomorphic to the 2-term $L_\infty$-algebra. In this case, the transferred structure is isomorphic to the induced one as 2-term $L_\infty$-algebras.

\begin{prop}
	\label{proplie2algeisomdecomp}
	Let $(L,d,[\cdot,\cdot],J)$ be a 2-term $L_\infty$-algebra, let $L'=L_0'\oplus L_1'$ be a graded vector space, and let $\phi\colon L\to L'$ be a graded linear isomorphism and let $\Phi\colon L_0\wedge L_0 \to L_1'$ be a linear map. Then $(L',d',[\cdot,\cdot]',J')$, with $d',[\cdot,\cdot]'$ and $J'$ uniquely defined by Equation \eqref{lowdimsemistrictlie2aglmoreq1}-\eqref{lowdimsemistrictlie2algmoreq4}, is a 2-term $L_\infty$-algebra isomorphic to $(L,d,[\cdot,\cdot], J)$. 
\end{prop}
\begin{myproof}
	Applying the Homotopy Transfer Theorem from e.g. \cite{loday2012algebraic} to the homotopy retract
	\begin{center}	
		\begin{tikzcd}
		{(L_0\oplus L_1,d)} \arrow["0",loop left, looseness=3] \arrow[r, "\phi", shift left] & {(L_0'\oplus L_1', d')} \arrow[l, "\phi^{-1}", shift left]
		\end{tikzcd}
	\end{center}
	yields that $(L',d',[\cdot,\cdot]',J')$ is a 2-term $L_\infty$-algebra and that $(\phi,\Phi)\colon L \to L'$ is a morphism of  2-term $L_\infty$-algebras. As $\phi$ is already an isomorphism of vector spaces, $(\phi,\Phi)$ is an isomorphism of 2-term $L_\infty$-algebras.
\end{myproof}
This proposition supports the intuition that transporting a given 2-term $L_\infty$-algebra via a supposed isomorphism of 2-term $L_\infty$-algebras induces another 2-term $L_\infty$-algebra, which is automatically isomorphic to the given one. Also note that for any linear graded isomorphism $\phi\colon L \to L'$, any linear map $\Phi\colon L_0 \wedge L_0 \to L_1'$ induces a 2-term $L_\infty$-algebra structure on $L'$ which are all isomorphic to $L$, as the maps $d'$, $[\cdot,\cdot]'$ and $J'$ will change accordingly.

\begin{theorem}
	\label{thmisomdecomp}
	Let $(L,d,[\cdot,\cdot],J)$ 2-term $L_\infty$-algebra with vector space decomposition $L=(\mathfrak{g}\oplus \im(d))\oplus (\ker(d)\oplus U)$. Then $(L,d,[\cdot,\cdot],J)$ is isomorphic to $L^{\mathfrak{g},\im(d),\rho,\tilde{J}}$.
\end{theorem}
\begin{myproof}
	Set $L':=L_0\oplus ( \ker(d)\oplus \im(d))$ and define $f$ and $h$ as in Lemma \ref{lem2exlie2alg}.
	By Proposition \ref{proplie2algeisomdecomp}, we find that the map
	\begin{align*}
	\begin{split}
	\phi:= \id \oplus f\colon L_0\oplus L_1 \to L_0' \oplus L_1'
	\end{split}
	\end{align*}
	is a graded linear isomorphism and that 
	\begin{align*}
	\begin{split}
	\Phi\colon L_0 \wedge L_0 &\to L_1'\\
	\Phi(x\wedge y)&:= ([x\supe{\im(d)},h(y)]+[x\supg,h(y)]+[h(x),y\supg])\supe{\ker(d)} +[x,y]\supe{\im(d)}
	\end{split}
	\end{align*}
	is a linear map.
	Note that the first term in this expression is antisymmetric by Equation \eqref{lowdimsemistrictdeflie2algeq1b}, as is the sum of the next two. 
	
	By Proposition \ref{proplie2algeisomdecomp}, we find that $(L',d',[\cdot,\cdot]',J')$, with $d',[\cdot,\cdot]'$ and $J'$ uniquely defined by Equation \eqref{lowdimsemistrictlie2aglmoreq1}-\eqref{lowdimsemistrictlie2algmoreq4}, is a 2-term $L_\infty$-algebra isomorphic to $(L,d,[\cdot,\cdot], J)$. We will now show that this 2-term $L_\infty$-algebra is identically $L^{\mathfrak{g},\im(d),\rho,\tilde{J}}$.
	
	For $x,y,z\in L_0$, $v\in L_1$ and $w\in L_1'$, we find that
	\begin{align*}
	[x,y]'= [\phi(x),\phi(y)] = \phi([x,y])-d'(\Phi (x\wedge y))= [x,y] - [x,y]\supe{\im(d)} = [x,y]\supg= [x\supg,y\supg]\supg,
	\end{align*}
	where the last equality follows from Equation \eqref{lowdimsemistrictdeflie2algeq1}.
	
	Moreover, we also obtain from Equation \eqref{lowdimsemistrictdeflie2algeq1} and Equation \eqref{lowdimsemistrictdeflie2algeq1b} that
	\begin{align*}
\begin{split}
	[f(v),y]'&= [\phi(v),\phi(y)]' = \phi([v,y])-\Phi (d(v)\wedge y)\\
	&= f([v,y]) - ([d(v)\supe{\im(d)},h(y)]+[d(v)\supg,h(y)]+[h(d(v)),y\supg])\supe{\ker(d)}-[d(v),y]\supe{\im(d)} \\
	&= [v,y]\supe{\ker(d)}+d([v,y]) - ([d(v),h(y)]+[h(d(v)),y\supg])\supe{\ker(d)} -[d(v),y]\\
	&= [v,y]\supe{\ker(d)} - ([h(d(v)),y\supe{\im(d)}]+[h(d(v)),y\supg])\supe{\ker(d)}\\
	&=[v-h(d(v)),y]\supe{\ker(d)}= [v\supe{\ker(d)},y]\supe{\ker(d)} = [v\supe{\ker(d)},y\supg]\supe{\ker(d)}.
\end{split}
	\end{align*}
	and thus that $[w,y]'= [w\supe{\ker(d)},y\supg]\supe{\ker(d)}$. Using those equations, we also obtain that
\begin{align*}
\begin{split}
\im(d) \ni dJ_{x,y,z}' = [x\supg,[y\supg,z\supg]\supg ]\supg &-[[x\supg,y\supg]\supg,z\supg]\supg - [y\supg,[x\supg,z\supg]\supg]\supg \in \mathfrak{g},
\end{split}
\end{align*}
	and as $\im(d) \cap \mathfrak{g} = 0$, we find that $dJ_{x,y,z}' =0$.
Similarly, we have that 
\begin{align*}
\begin{split}
J_{d(w),y,z}' &= [w\supe{\ker(d)},[y\supg,z\supg]\supg]\supe{\ker(d)} - [[w\supe{\ker(d)},y\supg]\supe{\ker(d)},z\supg]\supe{\ker(d)} - [y\supg,[w\supe{\ker(d)},z\supg]\supe{\ker(d)}]\supe{\ker(d)}\\
&= J_{d(w^{\ker(d)}),y,z}' =0.
\end{split}
\end{align*}
Combining the results from those two equations with Equation \eqref{lowdimsemistrictlie2algmoreq4}, we also find that
\begin{align*}
\begin{split}
J'_{x_1,x_2,x_3} &= J'_{x_1\supg,x_2\supg,x_3\supg}\supe{\ker(d)}\\
&= 
\Big(\phi(J_{x_1\supg,x_2\supg,x_3\supg})-
 \osum{\sigma \in \sh(1,2)}(-1)^\sigma[\phi(x_{\sigma(1)}\supg),\Phi_{x_{\sigma(2)}\supg,x_{\sigma(3)}\supg}]'+\Phi_{x_{\sigma(1)}\supg,[x_{\sigma(2)}\supg,x_{\sigma(3)}\supg]}\Big)\supe{\ker(d)}\\
&= J_{x_1\supg,x_2\supg,x_3\supg}\supe{\ker(d)}-\osum{\sigma \in \sh(1,2)}(-1)^\sigma [x_{\sigma(1)}\supg,h([x_{\sigma(2)}\supg,x_{\sigma(3)}\supg])]^{\ker(d)}.
\end{split}
\end{align*}
As the expressions for the differential, bracket and Jacobiator coincide, we conclude that $L'$ is identically $L^{\mathfrak{g},\im(d),\rho,\tilde{J}}$.
\end{myproof}

\begin{rem}
	\label{remskeletallie2algbaezshape}
	\textup{
		Given a general 2-term $L_\infty$-algebra $(L,d,[\cdot,\cdot],J)$, it is clear from both constructions that the skeletal 2-term $L_\infty$-algebra obtained from $(L,d,[\cdot,\cdot],J)$ in \cite{Baez2004}*{Proposition 51} is equal to $L^{\mathfrak{g},0,\rho,\tilde{J}}$, with $\mathfrak{g}$, $\rho$ and $\tilde{J}$ as in Lemma \ref{lem2exlie2alg}.
	}
\end{rem}
	We note that the above construction of $L^{\mathfrak{g},0,\rho,\tilde{J}}$ depends on a vector space decomposition of $L$. At first sight, this might lead to different 2-term $L_\infty$-algebras, but as Theorem \ref{lowdimthmlie2algdneq0liealgsnotisom} will show, this is not the case.
	
	In the next theorem, we will prove that two 2-term $L_\infty$-algebras in the shape constructed in Lemma \ref{lowdimlemexamplelie2alg} are isomorphic if and only if their underlying Lie algebras, vector spaces and representations are isomorphic, and if their 3-cocycles are cohomologous. Whereas it is quite straightforward to give a 2-term $L_\infty$-algebra isomorphism using the isomorphisms of the corresponding Lie algebras et cetera, it is a more intricate exercise to retrieve such morphisms form a 2-term $L_\infty$-algebra isomorphism. This is mostly due to the fact that a 2-term $L_\infty$-algebra isomorphism does not have to transport the Lie algebra underlying the first 2-term $L_\infty$-algebra to the Lie algebra underlying the second. We solve this problem by first considering the image of the Lie algebra under the 2-term $L_\infty$-algebra isomorphism and subsequently projecting onto the Lie algebra underlying the second 2-term $L_\infty$-algebra. 
\begin{theorem}
	\label{lowdimthmlie2algdneq0liealgsnotisom}
	Let $L:=L^{\mathfrak{g},U,\rho,\tilde{J}}$ and $L':=L^{\mathfrak{g'},U',\rho',\tilde{J'}}$ be two 2-term $L_\infty$-algebras from Lemma \ref{lowdimlemexamplelie2alg}. Then $L$ and $L'$ are isomorphic if and only if the following four statements hold:
	\begin{enumerate}
		\item $\mathfrak{g}$ and $\mathfrak{g'}$ are isomorphic Lie algebras;
		\item $U$ and $U'$ are isomorphic vector spaces;
		\item $(\rho,V)$ and $(\rho',V')$ are isomorphic representations;
		\item $\tilde{J}$ and $\tilde{J'}$ are cohomologous under the isomorphisms of 1. and 3.
	\end{enumerate}		
\end{theorem}
\begin{myproof}
	 In this proof, we denote the differential, bracket and Jacobiator of $L$ and of $L'$ by $d$, $[\cdot,\cdot]$, $J$ and by $d'$, $[\cdot,\cdot]'$, $J'$ respectively.
	 	 
$\Leftarrow$ First assume that 1.-4. hold. Then we have a Lie algebra isomorphism $\chi\colon \mathfrak{g}\to \mathfrak{g'}$, a linear isomorphism $f\colon U\to U'$, an intertwiner $g\colon V \to V'$ and a linear map $\tilde{\Phi}\colon\mathfrak{g}\wedge \mathfrak{g}\to V'$ such that
\begin{align*}
\begin{split}
g&(\tilde{J}_{x_1,x_2,x_3})-\tilde{J'}_{\chi(x_1),\chi(x_2),\chi(x_3)} = (\delta \tilde{\Phi})(x_1,x_2,x_3)
\end{split}
\end{align*}
for $x_1,x_2,x_3\in \mathfrak{g}$.\\
Then we define $\phi:= (\chi \oplus f) \oplus (g\oplus f) \colon L \to L'$ and
\begin{align*}
\begin{split}
\Phi\colon L_0\wedge L_0&\to L_1'\\
\Phi(x,y):&= \tilde{\Phi}(x\supg,y\supg).
\end{split}
\end{align*}
For $x,y,z \in L_0$ and $v\in L_1$, we find, by decoding the definitions of the given objects and maps, that
\begin{align*}
	\begin{split}
		\phi&([y,z])-[\phi(y),\phi(z)]'
		=\phi([y\supg,z\supg]_\mathfrak{g})-[\phi(y)\supe{\mathfrak{g'}},\phi(z)\supe{\mathfrak{g'}}]_\mathfrak{g'}\\
		&=\chi([y\supg,z\supg]_\mathfrak{g})-[\chi(y\supg),\chi(z\supg)]_\mathfrak{g'}=0 = d'(\tilde{\Phi}_{y\supg,z\supg}) = d'(\Phi_{y,z}).
	\end{split}
\end{align*}
and that
\begin{align*}
	\begin{split}
		\phi([v,y]) &= -\phi(\rho(y\supg)v\supv)
		=-g(\rho(y\supg)v\supv)\\
		&=-\rho'(\chi(y\supg))g(v\supv)
		=-\rho'(\phi(y)\supe{\mathfrak{g'}})\phi(v)\supe{V'}
		=[\phi(v),\phi(y)]'.
	\end{split}
\end{align*}
The latter equation implies that
\begin{align*}
\begin{split}
	\phi([v,y])-[\phi(v),\phi(y)]' = 0 = \tilde{\Phi}_{d(v)\supg,y\supg} = \Phi_{d(v),y}.
\end{split}
\end{align*}
Moreover, it is not hard to see that further unfolding of definitions yields that
\begin{align*}
	\tilde{\Phi}_{[x\supg,y\supg]_\mathfrak{g},z\supg} = \tilde{\Phi}_{[x,y]\supg,z\supg} = \Phi_{[x,y],z}
\end{align*}
and that
\begin{align*}
	\begin{split}
		\rho'(\phi(x\supg))\tilde{\Phi}_{y\supg,z\supg}
		=\rho'(\phi(x)\supe{\mathfrak{g'}})\tilde{\Phi}_{y\supg,z\supg}
		=[\phi(x),\tilde{\Phi}_{y\supg,z\supg}]'
		= [\phi(x),\Phi_{y,z}]'.
	\end{split}
\end{align*}
Combining those last two results gives us that
\begin{align*}
	\begin{split}
		\phi&(J_{x_1,x_2,x_3}) - J'_{\phi(x_1),\phi(x_2),\phi(x_3)}
		=g(\tilde{J}_{x_1\supg,x_2\supg,x_3\supg})- \tilde{J'}_{\phi(x_1)^{\mathfrak{g'}} ,\phi(x_2)^{\mathfrak{g'}},\phi(x_3)^{\mathfrak{g'}}}\\
		&=g(\tilde{J}_{x_1\supg,x_2\supg,x_3\supg})- \tilde{J'}_{\chi(x_1\supg),\chi(x_2\supg),\chi(x_3\supg)}= (\delta \tilde{\Phi})(x_1\supg,x_2\supg,x_3\supg)\\
		&=\osum{\sigma \in \sh(1,2)}(-1)^\sigma\rho'(\phi(x_{\sigma(1)}\supg))\tilde{\Phi}_{x_{\sigma(2)}\supg,x_{\sigma(3)}\supg} -\osum{\sigma \in \sh(2,1)}(-1)^\sigma \tilde{\Phi}_{[x_{\sigma(1)}\supg,x_{\sigma(2)}\supg]_\mathfrak{g},x_{\sigma(3)}\supg}\\		
		&= \osum{\sigma \in \sh(1,2)}(-1)^\sigma\big([\phi(x_{\sigma(1)}),\Phi_{x_{\sigma(2)},x_{\sigma(3)}}]'+\Phi_{x_{\sigma(1)},[x_{\sigma(2)},x_{\sigma(3)}]}\big)
	\end{split}
\end{align*}
for $x_1,x_2,x_3\in L_0$. The above calculations prove that $(\phi,\Phi)$ is a 2-term $L_\infty$-algebra morphism, and as $\phi$ is invertible, it is even an isomorphism of 2-term $L_\infty$-algebras. 

$\Rightarrow$ Now assume that $L$ and $L'$ are isomorphic 2-term $L_\infty$-algebras, and let $(\phi,\Phi)$ denote an isomorphism between them. Note that by Equation \eqref{lowdimsemistrictlie2aglmoreq1}, we have for $u\in U$ in degree 0 that
\begin{align*}
	\phi(u) = \phi(d(u))=d'(\phi(u)),
\end{align*}
which proves that $\phi(U)\subseteq U'$ in degree 0. Analogously, we find that $\phi^{-1}(U')\subseteq U$, which proves that $\phi(U)= U'$ and thus that $U$ and $U'$ are isomorphic vector spaces. 

If we let $v\in V$, we obtain that $d'(\phi(v)) = \phi(d(v))=\phi(0)=0$, which implies that $\phi(v) \in V'$. Analogously, we find that $\phi^{-1}(v')\in V$ for $v'\in V'$, and combining this with the previous observation gives us that $\phi(V)=V'$.

As this implies that $\mathfrak{g}'\oplus U' = L_0' = \phi(\mathfrak{g}\oplus U) = \phi(\mathfrak{g})\oplus U'$, it is a simple exercise in linear algebra to deduce that
\begin{align*}
\begin{split}
\pi\colon \phi(\mathfrak{g}) &\to \mathfrak{g}'\\
x &\mapsto x\supe{\mathfrak{g'}} \quad \text{(where we view\ } x\in \mathfrak{g}'\oplus U' \text{)}
\end{split}
\end{align*}
is a linear isomorphism. If we let $x,y\in \mathfrak{g}$, we have that $\phi([x,y])-\pi(\phi([x,y]))\in U'$, and thus that
\begin{align}
	\label{eqthmisomdphi}
	\begin{split}
			U'\ni d'(\Phi(x,y)) &=\phi([x,y])  -[\phi(x),\phi(y)]'\\
			&= 	\phi([x,y])  -[\phi(x),\phi(y)]'+\pi(\phi([x,y]))-\pi(\phi([x,y]))\\
		&=\underbrace{\pi(\phi([x,y]))-[\phi(x),\phi(y)]'}_{\in \mathfrak{g'}} \oplus  \underbrace{\phi([x,y]) -\pi(\phi([x,y]))}_{\in U'}.
	\end{split}
\end{align}
By the uniqueness of the direct sum decomposition, we obtain that
\begin{align*}
	\begin{split}
		\pi(\phi([x,y]))-[\phi(x),\phi(y)]'=0.
	\end{split}
\end{align*}
This gives us that
\begin{align*}
	\begin{split}
		\pi(\phi([x,y]_{\mathfrak{g}})) &=\pi(\phi([x,y])) = [\phi(x),\phi(y)]' = [\phi(x)^{\mathfrak{g'}},\phi(y)^{\mathfrak{g'}}]_{\mathfrak{g'}}
		= [\pi(\phi(x)),\pi(\phi(y))]_{\mathfrak{g'}}.
	\end{split}
\end{align*}
We thus obtain that $\mathfrak{g}$ and $\mathfrak{g'}$ are isomorphic Lie algebras, as 
\begin{align*}
	\tau:=\pi \circ \phi|_{\mathfrak{g}}\colon \mathfrak{g} \to \mathfrak{g'}
\end{align*}
is a Lie algebra isomorphism. If we let $v\in V$ and $y \in \mathfrak{g}$, we also find that
\begin{align*}
	\begin{split}
		0=\Phi(d(v),y)=\phi([v,y])-[\phi(v),\phi(y)]' = -\phi(\rho(y)v) + \rho'(\phi(y)^{\mathfrak{g'}})\phi(v),
	\end{split}
\end{align*}
that is,
\begin{align*}
	\begin{split}
	\phi(\rho(y)v) = \rho'(\phi(y)^{\mathfrak{g'}})\phi(v)
	= \rho'(\pi(\phi(y)))\phi(v)
		=\rho'(\tau(y))\phi(v).
	\end{split}
\end{align*}
This proves that $\rho$ and $\rho'$ are isomorphic representations.

Furthermore, for $x_1,x_2,x_3 \in \mathfrak{g}$, we have that
\begin{align*}
	[\phi(x_1),\Phi_{x_2,x_3}]' = \rho'(\phi(x_1)^{\mathfrak{g'}})\Phi_{x_2,x_3}\supe{V'} =\rho'(\tau(x_1))\Phi_{x_2,x_3}\supe{V'}.
\end{align*}
By Equation \eqref{eqthmisomdphi} and the Jacobi identity on $\mathfrak{g}$, we obtain that
\begin{align*}
	\begin{split}
		d'\big(\osum{\sigma \in \sh(2,1)}(-1)^\sigma \Phi_{[x_{\sigma(1)},x_{\sigma(2)}],x_{\sigma(3)}}\big)
		&= \osum{\sigma \in \sh(2,1)}(-1)^\sigma \phi([[x_{\sigma(1)},x_{\sigma(2)}],x_{\sigma(3)}])\supe{U'}\\
		&= \phi\big(\osum{\sigma \in \sh(2,1)}(-1)^\sigma [[x_{\sigma(1)},x_{\sigma(2)}]\subg,x_{\sigma(3)}]\subg\big)\supe{U'} =0,
	\end{split}
\end{align*}
from which we conclude that
\begin{align*}
	\begin{split}
		\phi(\tilde{J}_{x_1,x_2,x_3}) &- \tilde{J}'_{\phi(x_1),\phi(x_2),\phi(x_3)}\\
		&= \osum{\sigma \in \sh(1,2)}(-1)^\sigma[\phi(x_{\sigma(1)}),\Phi_{x_{\sigma(2)},x_{\sigma(3)}}]'-\osum{\sigma \in \sh(2,1)}(-1)^\sigma \Phi_{[x_{\sigma(1)},x_{\sigma(2)}],x_{\sigma(3)}}\\
		&= 	\osum{\sigma \in \sh(1,2)}(-1)^\sigma\rho'(\tau(x_{\sigma(1)}))\Phi_{x_{\sigma(2)},x_{\sigma(3)}}\supe{V'}
-\osum{\sigma \in \sh(2,1)}(-1)^\sigma \Phi_{[x_{\sigma(1)},x_{\sigma(2)}]_\mathfrak{g},x_{\sigma(3)}}\supe{V'}\\
&= (\delta (r\circ \Phi|_{\mathfrak{g}\wedge \mathfrak{g}}))(x_1,x_2,x_3),
	\end{split}	
\end{align*}
where $r\colon L_1'\to V'$ is the projection. This proves that $\tilde{J}$ and $\tilde{J}'$ are cohomologous.

\end{myproof}
For a given 2-term $L_\infty$-algebra, \cite{Baez2004}*{Proposition 51} gives an equivalent skeletal 2-term $L_\infty$-algebra and \cite{Baez2004}*{Theorem 55} yields a Lie algebra $\mathfrak{g}$, Lie algebra representation $\rho$ of $\mathfrak{g}$ and 3-cocycle of $\rho$ on $\mathfrak{g}$ corresponding to that skeletal 2-term $L_\infty$-algebra. Similarly, Theorem \ref{thmisomdecomp}  also yields a Lie algebra, Lie algebra representation and 3-cocycle corresponding to the same 2-term $L_\infty$-algebra. By Remark \ref{remskeletallie2algbaezshape}, those Lie algebras, Lie algebra representations and 3-cocycles coincide. As \cite{Baez2004}*{Theorem 57} states that two skeletal 2-term $L_\infty$-algebras are equivalent if and only if their corresponding Lie algebras and Lie algebra representations are isomorphic and their corresponding 3-cocycles are cohomologous, we can use Theorem \ref{lowdimthmlie2algdneq0liealgsnotisom} to come to the following conclusion.
\begin{cor}
	\label{corequivisomlie2alg}
	Two equivalent 2-term $L_\infty$-algebras are isomorphic if and only if the images of their differentials are isomorphic vector spaces.\qed
\end{cor}
Using Corollary \ref{corequivisomlie2alg}, we can immediately conclude that the skeletal string 2-term $L_\infty$-algebra from Example \ref{exampleskeletalstringlie2alg} and the strict string 2-term $L_\infty$-algebra from Example \ref{examplestrictstringlie2alg} are not isomorphic, as the image of the differential of the former is zero-dimensional, and the image of the latter has an uncountably many basis elements. Another equivalent 2-term $L_\infty$-algebra can be found in \cite{Wagemann2006}, and it too can be shown to be non-isomorphic to either of the previous examples by looking at the image of the respective differential.

Determining if 2-term $L_\infty$-algebras are isomorphic can by done by either using \cite{Baez2004}*{Theorem 57} and Corollary \ref{corequivisomlie2alg}, or by Lemma \ref{lem2exlie2alg}, Theorem \ref{thmisomdecomp} and Theorem \ref{lowdimthmlie2algdneq0liealgsnotisom}. This completes the classification procedure. 

\section{Discussion}
\label{secdis}
In this article, we have classified all 2-term $L_\infty$-algebras, and thus also all semistrict Lie 2-algebras \cite{Baez2004}. We have found that a 2-term $L_\infty$-algebra is a combination of a Lie algebra with a vector space, a representation and a cohomology class. This implies that all further classification for 2-term $L_\infty$-algebras can be done solely in terms of Lie algebras and their representation theory. The procedure explained in this article can be interpreted as an extension of a former approach to obtain a clear view of all 2-term $L_\infty$-algebras in terms of more familiar and more studied objects. Even though the isomorphism classes of 2-term $L_\infty$-algebras are more restrictive than the equivalence classes, the only additional data that is required to distinguish isomorphism classes is the vector space that is modded out by taking homology. Hence, it can be concluded that the category of 2-term $L_\infty$-algebras is an enrichment of the category of Lie algebras, and this can simplify our thinking about 2-term $L_\infty$-algebras and thus about semistrict Lie 2-algebras as well. 

The question remains, though, if it is possible to find a similar type of classification for hemistrict Lie 2-algebras or more general $L_\infty$-algebras. In a further attempt to classify more $L_\infty$-algebras, the suggested direction would be to classify all 3-term $L_\infty$-algebras in a similar vein as the classification in this article. This is more involved, as a 3-term $L_\infty$-algebra has a (non-zero) differential in two different degrees, so, by taking homology, there are multiple vector spaces that are modded out which have to be taken into account. Moreover, morphisms of 3-term $L_\infty$-algebras contain a new map of degree 2 which has to satisfy more equalities and which does at least at first sight not appear to allow for a characterization of 3-term $L_\infty$-algebras as simple as the characterization of 2-term $L_\infty$-algebras.

This article could still shed a light on possible ways creating new examples and classifications of more general $L_\infty$-algebras and we hope that this article assists in further understanding and classification of $L_\infty$-algebras in the near future.

\paragraph{Acknowledgements}
The author is funded by the Fundamentals of the Universe program at the University of Groningen. The author would like to thank Roland van der Veen and Jorge Becerra for fruitful discussions and for their assistance with the final revisions.

\bibliographystyle{amsplain}

\end{document}